\documentclass[10pt]{article}

\DeclareMathAlphabet{\mathpzc}{OT1}{pzc}{m}{it}

\usepackage{amsmath}
\usepackage{amssymb}
\usepackage{amsfonts}
\usepackage{amsthm}
\usepackage{mathrsfs}
\usepackage{ifsym}
\usepackage{fontenc}
\usepackage{textcomp}
\usepackage{tipa}
\usepackage{enumerate}
\usepackage[pdftex]{color, graphicx}

\begin{document}
 
\title{{\bf  Supersingular conjectures for the Fricke group}}       
\author{Patrick Morton}        
\date{May 30, 2021}          
\maketitle

\begin{abstract} A proof is given of several conjectures from a recent paper of Nakaya concerning the supersingular polynomial $ss_p^{(N*)}(X)$ for the Fricke group $\Gamma_0^*(N)$, for $N \in \{2, 3, 5, 7\}$.  One of these conjectures gives a formula for the square of $ss_p^{(N*)}(X)$ (mod $p$) in terms of a certain resultant, and the other relates the primes $p$ for which $ss_p^{(N*)}(X)$ splits into linear factors (mod $p$) to the orders of certain sporadic simple groups.
\end{abstract}

\section{Introduction.}

This paper is devoted to proving several of the conjectures appearing in Nakaya's paper \cite{na}.  These conjectures concern the supersingular polynomial $ss_p^{(N*)}(X)$ for the Fricke group $\Gamma_0^*(N)$, where $N \in \{2, 3, 5, 7\}$.  This polynomial is defined as follows.  (See \cite[p. 2254]{sa2} and \cite[p. 486]{na}.)  Define the polynomials
\begin{align*}
R_2(X,Y) &= X^2-X(Y^2-207Y+3456)+(Y+144)^3,\\
R_3(X,Y) &= X^2 - XY(Y^2-126Y+2944)+Y(Y+192)^3,\\
R_5(X,Y) &= X^2-X(Y^5-80Y^4+1890Y^3-12600Y^2+7776Y+3456)\\
& \ \ +(Y^2+216Y+144)^3,\\
R_7(X,Y) &= X^2-XY(Y^2-21Y+8)(Y^4-42Y^3+454Y^2-1008Y-1280)\\
& \ \ +Y^2(Y^2+224Y+448)^3.
\end{align*}
For each $N$ and each prime $p \neq N$, the polynomial $R_N(X,Y)$ is absolutely irreducible over $\mathbb{F}_p$ and defines a correspondence (in the sense of algebraic geometry) between the projective line $\mathbb{P}^1(\overline{\mathbb{F}}_p)$ and itself; or equivalently, between the rational function field $\overline{\mathbb{F}}_p(x)$ and itself (see \cite{de}).  In this correspondence, the points in $\overline{\mathbb{F}}_p$ corresponding to the $j$-invariants of supersingular elliptic curves are the supersingular invariants $j^*$ for $\Gamma_0^*(N)$, and they are roots of a polynomial $ss_p^{(N*)}(X) \in \mathbb{F}_p[X]$:
$$ss_p^{(N*)}(X) = \prod_{R_N(j,j^*)=0, ss_p(j)=0}{(X-j^*)} \in \mathbb{F}_p[X];$$
the product running over the distinct roots $j^*$ of $R_N(j,j^*) = 0$ in $\overline{\mathbb{F}}_p$, as $j$ runs over the supersingular $j$-invariants in characteristic $p$.  (See \cite{sa1}, \cite{sa2}, \cite{na}.)  It is well-known that the roots of $ss_p(X)$ lie in $\mathbb{F}_{p^2}$, and it was shown in \cite[Section 6]{mor1} that the values $j^*$ lie in $\mathbb{F}_{p^2}$, for $N \in \{2, 3, 5, 7\}$.  (See Tables 2 and 3 below for $N=5, 7$.)  Consequently, the above correspondence can be defined over $\mathbb{F}_{p^2}$.  \medskip

Nakaya's Conjecture 4 takes the general form
$$A_{N,p}(Y) \textrm{Res}_X(R_N(X,Y), ss_p(X)) = B_{N,p}(Y) ss_p^{(N*)}(Y)^2 \ \ (\textrm{mod} \ p),$$
where $A_{N,p}(Y)$ and $B_{N,p}(Y)$ are explicit polynomials of low degree which depend on $N$ and the residue class of $p$ (mod $12N$).  This formula arises from the fact that the correspondence $X \rightarrow Y$ is generally $2-1$, i.e. two values of $X$ correspond to a single value of $Y$.  Exceptions occur where the correspondence is ``ramified", i.e., when $j = 0$ or $j=1728$ is supersingular in characteristic $p$, and for several other values of $j$ in characteristic $p$, depending on $N$. \medskip

The proof of the above formula requires knowing a number of ring class polynomials $H_{d}(X)$ explicitly (see \cite{co}), and uses Deuring's fundamental theorem \cite{de1} that if  $\mathfrak{j}$ is the $j$-invariant of an elliptic curve in characteristic $0$ with complex multiplication by the imaginary quadratic order $\mathcal{O} = \textsf{R}_{d}$ of discriminant $d$, then the reduction $j \equiv \mathfrak{j}$ (\textrm{mod} $\mathfrak{p}$), modulo a prime divisor $\mathfrak{p}$ of $p$, is supersingular if and only if the Legendre symbol $\left(\frac{d}{p}\right) \neq 1$.  Thus, part of the proof involves recognizing several ring class polynomials and their associated discriminants.  See Lemmas 2 and 4 and their proofs.  The proof also requires the fact that two isogenous elliptic curves in characteristic $p$ are both supersingular when one of them is (see \cite{si}).  In the cases $N=5, 7$ this requires that we exhibit an explicit isogeny between the Tate normal form $E_N$ for a point of order $N$ and its isogenous curve $E_{N,N} =E_N/\langle (0,0) \rangle$, in order to calculate their $j$-invariants.  For $N=5$ this isogeny has been worked out in \cite{mor0} and \cite{mor}.  For $N=7$ we give a summary of the calculation in Section 3 (see Fact 7). \medskip

In Section 2 (Theorem 1) we work out the case $N=5$, and in Section 3 (Theorem 3) we deal with the case $N=7$.  The cases $N=2,3$ are handled in Section 4 (Theorem 5).  Taken together, these theorems cover all four cases of Nakaya's Conjecture 4.  \medskip

In Section 5 we give a simple proof of Nakaya's Conjecture 2 \cite{na}, which says that in the cases $N=5, 7$ the primes for which $ss_p^{(N*)}(X)$ splits into linear factors over $\mathbb{F}_p$ coincide with the prime divisors of the order of a specific sporadic simple group (the Harada-Norton group $HN$ and the Held group $He$, respectively; see \cite[Ch. 10]{con}).  Nakaya proved the analogous results for $N=2,3$ in \cite{na} using an explicit formula for the number of linear factors of $ss_p^{(N*)}(X)$ over $\mathbb{F}_p$ and a class number estimate.  The proof given in Theorem 6 below is elementary, uses nothing about class numbers, and is also valid for the cases $N=2,3$ discussed by Nakaya.  It shows that the set of primes for which $ss_p^{(N*)}(X)$ splits into linear factors modulo $p$ is always a subset of the primes for which the supersingular polynomial $ss_p(X)$ splits (mod $p$), so that the proof of Nakaya's Conjecture 2 requires only a modest calculation.

\section{The case $N=5$.}

Let the curve $R_5(X,Y) = 0$ be defined by
\begin{align*}
R_5(X,Y) & = X^2-X(Y^5-80Y^4+1890Y^3-12600Y^2+7776Y+3456)\\
& \ \ +(Y^2+216Y+144)^3.
\end{align*}
This is a curve of genus $0$ parametrized by the equations
$$X = -\frac{(z^2+12z+16)^3}{z+11}, \ \ \ Y = -\frac{z^2+4}{z+11}.$$
See \cite[p. 263]{mor1}.  We have
\begin{align}
\textrm{disc}_X R_5(X,Y) &= Y^2(Y-4)^2(Y-18)^2(Y-36)^2(Y^2-44Y-16)\\
\textrm{disc}_Y R_5(X,Y) & = 5^5 X^4 (X-1728)^4 (X+32^3)^2(X-66^3)^2(X+96^3)^2\\
\notag & = 5^5X^4 (X-1728)^4 H_{-11}(X)^2 H_{-16}(X)^2 H_{-19}(X)^2.
\end{align}
Define
\begin{align}
\mu_5 &= \frac{1}{2}\left(1-\left(\frac{-5}{p}\right)\right)\\
\delta & = \frac{1}{2}\left(1-\left(\frac{-3}{p}\right)\right)\\
\varepsilon & = \frac{1}{2}\left(1-\left(\frac{-4}{p}\right)\right).
\end{align}

In this section we will prove the following theorem, conjectured by Nakaya \cite[Conjecture 4]{na}. \bigskip

\noindent {\bf Theorem 1.} {\it If $p \ge 7$ is a prime and $ss_p(X)$ denotes the supersingular polynomial in characteristic $p$, then}
\begin{align}
(Y^2-&44Y-16)^{\mu_5} \textrm{Res}_X(ss_p(X),R_5(X,Y)) \equiv \\
& (Y^2+216Y+144)^{2\delta}(Y^2-540Y-6480)^\varepsilon ss_p^{(5*)}(Y)^2 \ (\textrm{mod} \ p).
\end{align}
\medskip

\noindent {\bf Lemma 2.} {\it We have the following class equations:}
\begin{align*}
H_{-20}(X) &= X^2-1264000X-681472000,\\
H_{-75}(X) & = X^2 +654403829760X+5209253090426880,\\
H_{-100}(X) & = X^2-44031499226496X-292143758886942437376.
\end{align*}

\noindent {\it Proof.}. For $H_{-20}(X)$, see \cite[p. 8]{mor2}.  For $H_{-75}(X)$, note from \cite[p. 311]{ber} that
\begin{align*}
\gamma_2 \left(\frac{3+\sqrt{-75}}{2}\right) &= \left(j\left(\frac{3+\sqrt{-75}}{2}\right)\right)^{1/3}\\
&=-32 \cdot 3 \cdot 5^{1/6} \left(\frac{69+31\sqrt{5}}{2}\right).
\end{align*}
Then $H_{-75}(X)$ is the minimal polynomial of the quadratic irrationality
$$j\left(\frac{3+\sqrt{-75}}{2}\right) = -32^3 \cdot 3^3 \cdot \sqrt{5} \left(\frac{69+31\sqrt{5}}{2}\right)^3.$$
To compute the class equation $H_{-100}(X)$ we use the Rogers-Ramanujan continued fraction $r(\tau)$.  From \cite[p. 138]{du} we have the well-known value of Ramanujan,
$$r(i) = \sqrt{\frac{5+\sqrt{5}}{2}}-\frac{1+\sqrt{5}}{2},$$
whose minimal polynomial is
$$f(x) = x^4 + 2x^3 - 6x^2 - 2x + 1.$$
The minimal polynomial $f_5(x)$ of $r(5i)$ can be found using the identity
$$r^5(\tau) = r \frac{r^4-3r^3+4r^2-2r+1}{r^4+2r^3+4r^2+3r+1}, \ \ \ r = r(5\tau).$$
See \cite[p. 93]{anb}.  Setting
$$g(x,y) = (y^4+2y^3+4y^2+3y+1)x^5 - y(y^4-3y^3+4y^2-2y+1),$$
the polynomial $f_5(x)$ must divide the resultant
\begin{align*}
\textrm{Res}_t&(f(t),g(t,x)) = x^{20} + 510x^{19} - 13590x^{18} + 32280x^{17} - 82230x^{16} + 153522x^{15}\\
& \ \  - 302910x^{14} + 273540x^{13} - 412830x^{12} + 268230x^{11} - 262006x^{10}\\
& \ \ - 268230x^9 - 412830x^8- 273540x^7 - 302910x^6 - 153522x^5\\
& \ \ - 82230x^4 - 32280x^3 - 13590x^2 - 510x + 1.
\end{align*}
This polynomial is irreducible, and so equals $f_5(x)$.  Now $j(5i)$ can be found from the relation
$$F(r,j) = (r^{20}-228r^{15}+494r^{10}+228r^5+1)^3+j(\tau) r^5 (r^{10}+11r^5-1)^5 = 0.$$
(See \cite[p. 138]{du}.). Taking the resultant
$$\textrm{Res}_t(f_5(t),F(t,X)) = 5^{300} (X^2 - 44031499226496X - 292143758886942437376)^{10}$$
shows that $H_{-100}(X)$, which is the minimal polynomial of $j(5i)$, is given by the polynomial in the lemma.   
See also the values for $j(5i)$ given in \cite[p. 58]{berw} and \cite[p. 202]{fr1}.  \smallskip

A similar proof may be given for $H_{-75}(X)$ starting with the value $r(\rho)$ in place of $r(i)$, where $\rho = \frac{-1+\sqrt{3}i}{2}$: 
$$r(\rho) = e^{-\pi i/5} \left(\frac{\sqrt{30+6\sqrt{5}}-3-\sqrt{5}}{4}\right),$$
whose fifth power has the minimal polynomial
$$g_3(x) = x^4 - 228x^3 + 494x^2 + 228x + 1.$$
See \cite[Eq. (2.4)]{du}.  $\square$ \bigskip

The proof of Theorem 1 is given in the course of verifying the following facts. \medskip

\noindent {\bf Fact 1}. {\it All the roots of $(Y^2-44Y-16)^{\mu_5} \textrm{Res}_X(ss_p(X),R_5(X,Y))$ are roots of $ss_p^{(5*)}(Y)$.}

This is clear by definition for the resultant.  The factor $Y^2-44Y-16$ arises from roots $X$ of $H_{-20}(X)$, since
$$\textrm{Res}_Y(R_5(X,Y),Y^2-44Y-16) = (X^2-1264000X-681472000)^2 = H_{-20}(X)^2.$$
Furthermore,
\begin{equation}
\textrm{Res}_X(H_{-20}(X),R_5(X,Y))  = (Y^2-44Y-16) h_{20}(Y),
\end{equation}
where
\begin{align*}
h_{20}(Y) & =  Y^{10} - 1262660Y^9 - 1454280320Y^8 - 69099329600Y^7\\
& - 10276940953600Y^6 + 460141172243456Y^5 - 3888238950420480Y^4\\
&  - 12956776173404160Y^3 - 415176163957145600Y^2\\
& - 345243549014425600Y - 512182838955606016.
\end{align*}
Since the roots of $H_{-20}(X)$ are supersingular $j$-invariants in characteristic $p$ exactly when $\left(\frac{-5}{p}\right) = -1$ (for primes $p > 7$), i.e., when $\mu_5 = 1$, we see that the roots of $Y^2-44Y-16=0$ are roots of $ss_p^{(5*)}(X)$ in this case. \medskip

\noindent {\bf Fact 2.} {\it Since $R_5(X,Y)$ is quadratic in $X$, each root $y$ of $ss_p^{(5)}(Y) = 0$ arises from exactly two roots $x$ of $R_5(X,y)=0$, except for the following values, which are all roots of the discriminant in equation (1).}
\begin{align*}
y &= 0 \ \textrm{corresponds to} \ x = 1728, \ \textrm{since} \ R_5(1728,Y) = Y^2h_4(Y)^2\\
& =  Y^2(Y^2-540Y-6480)^2 \ \textrm{and} \ R_5(X,0) = (X-1728)^2;\\
y &= 4 \ \textrm{corresponds to} \ x = -32^3, \ \textrm{since} \ R_5(-32^3,Y) = (Y-4)^2 h_{11}(Y)\\
&  = (Y-4)^2(Y^4 + 33424Y^3- 2213664Y^2 + 53951744Y + 74373376)\\
&  \ \textrm{and} \ R_5(X,4) = (X+32^3)^2;\\
y &= 18 \ \textrm{corresponds to} \ x = 66^3, \ \textrm{since} \ R_5(66^3,Y) = (Y-18)^2 h_{16}(Y)\\
& =(Y-18)^2(Y^4 - 286812Y^3 +12814524 Y^2+21146832 Y+252047376)\\
&  \ \textrm{and} \ R_5(X,18) = (X-66^3)^2;\\
y &= 36 \ \textrm{corresponds to} \ x = -96^3, \ \textrm{since} \ R_5(-96^3,Y) = (Y-36)^2h_{19}(Y)\\
& = (Y-36)^2(Y^4 + 885456Y^3 -6886944Y^2 + 39004416Y + 606341376)\\
&  \ \textrm{and} \ R_5(X,36) = (X+96^3)^2.
\end{align*}
It follows that for these values $(Y-y)^2$ exactly divides the resultant in (6), when the corresponding $X$-value is supersingular (corresponding to quadratic discriminants $d = -4, -11, -16, -19$, see (2)), and so are accounted for in (7) by the factor $ss_p^{(5*)}(Y)^2$.  This requires that we take $p$ to be a prime not dividing the values at $y$ of each of the cofactors of $(Y-y)^2$ in these four cases:
\begin{align*}
y=0: \ \ &h_4(0) = -6480 = -2^4 \cdot 3^4 \cdot 5;\\
y= 4: \ \ & h_{11}(4) = 256901120 = 2^{20} \cdot 5 \cdot 7^2;\\
y=18: \ \ & h_{16}(18) = 3112013520 = 2^4 \cdot 3^8 \cdot 5 \cdot 7^2 \cdot 11^2;\\
y=36: \ \ & h_{19}(36) = 34398535680 = 2^{20} \cdot 3^8 \cdot 5.
\end{align*}
Hence, we must require that $p \notin \{2,3,5,7,11\}$. \medskip

Finally, each of the roots of $Y^2-44Y-16$ arises from only one value of $X$, by the first resultant calculation in Fact 1.  The second resultant calculation (8) shows that this factor occurs only to the first power in $\textrm{Res}_X(ss_p(X),R_5(X,Y))$, when $p \notin \{2,5,11,13,17,19\}$, which is the set of primes dividing the integer resultant
$$\textrm{Res}_Y(Y^2-44Y-16,h_{20}(Y)) = 2^{60} \cdot 5^6 \cdot 11^6 \cdot 13^4 \cdot 17^4 \cdot19^2.$$
Hence, $Y^2-44Y-16$ and $h_{20}(Y)$ have no factor in common when $p>19$ and $\mu_5 = 1$; then the factor $(Y^2-44Y-16)^2$ exactly divides (6) and is accounted for by the same factor of $ss_p^{(5*)}(Y)^2$ in (7).  Otherwise, $\mu_5=0$ and the roots of $H_{-20}(X)$ are not supersingular for $p$, in which case the factor $Y^2-44Y-16$ does not occur.  \medskip

Note that the $Y$-values above are distinct for primes $p > 19$, since
$$\textrm{disc}_Y(Y(Y - 4)(Y - 18)(Y - 36)(Y^2 - 44Y - 16)) = 2^{56} \cdot 3^{12} \cdot 5^3 \cdot 7^2 \cdot 11^6 \cdot 19^2.$$
Similarly, the above $X$-values, i.e. the roots of (2), are distinct for $p>19$ and $p \neq 43, 67$, since
\begin{align*}
\textrm{disc}_X(X(X - 1728)(X + 32^3)(X - 66^3)(X+96^3)) = & \ 2^{152} \cdot 3^{56} \cdot 7^{12} \cdot 11^8 \cdot 13^2\\
 & \cdot 19^4 \cdot 43^2 \cdot 67^2.
\end{align*}

\noindent {\bf Fact 3.} {\it If $y$ is a root of (6) corresponding to two distinct $X$-values, and one of these values $x$ is a root of $ss_p(X)$, then the second value is also.}

This can be seen as follows.  It suffices to show this for the resultant in (6).  It can be checked on Maple that
\begin{equation}
R_5\left(X,-\frac{z^2+4}{z+11}\right) = \left(X+\frac{(z^2+12z+16)^3}{z+11}\right) \left(X+\frac{(z^2-228z+496)^3}{(z+11)^5}\right).
\end{equation}
By \cite[Eqs. (5), (8)]{mor} and \cite[pp. 258-259]{mor0}, the roots of (9), namely
$$j_5 = -\frac{(z^2+12z+16)^3}{z+11} \ \ \textrm{and} \ \ j_{5,5} = -\frac{(z^2-228z+496)^3}{(z+11)^5},$$
are the $j$-invariants of the isogenous elliptic curves
$$E_5: \ Y^2+(1+b)XY+bY = X^3 +bX^2, \ \ \ z = b-\frac{1}{b},$$
(this is the Tate normal form for a point of order $5$) and
$$E_{5,5}: \ Y^2 +(1+b)XY+ 5bY = X^3+7bX^2+6(b^3+b^2-b)X+b^5+b^4-10b^3-29b^2-b.$$
If $j_5$ is supersingular, then $j_{5,5}$ is supersingular, and vice-versa.
\medskip

\noindent {\bf Fact 4.} {\it The only roots $y$ of $ss_p^{(5*)}(Y)$ which occur to a power higher than the second in (6) are those which correspond to the roots of the discriminant (2), because $R_5(x,Y)$ must have the square factor $(Y-y)^2$ for at least one of the $X$-values $x$ corresponding to $Y = y$.}

We have already discussed these roots in Fact 2, except for $x=0$.  We can igonore the cofactors $h_{11}(Y), h_{16}(Y), h_{19}(Y)$ in Fact 2 for $x = -32^3, 66^3, -96^3$ (corresponding to $y = 4, 18, 36$) when the prime $p \notin \{2, 3, 5, 7, 11, 13, 19, 43, 67\}$, since this set contains the prime factors of the discriminants of these polynomials (as well as the discriminants of $h_4(Y)$ and $Y^2+216Y+144$; see below).  For all other primes, these cofactors do not have multiple roots; and since their factors do not occur to a power higher than the first for the other roots $x \in \{-32^3, 66^3,-96^3\}$ of (2), they cannot occur to a power higher than the second in (6), unless one of these roots also occurs for $x=0$ or $x=1728$.  Any such roots will be covered by the cases $x=1728$ and $x=0$ considered next.

The multiple roots $y$ corresponding to $x=1728$ in Fact 2 come from the factorization
$$R_5(1728,Y) = Y^2(Y^2-540Y-6480)^2.$$
Note that
\begin{align*}
\textrm{Res}_Y(&R_5(X,Y), Y^2-540Y-6480) \\
& = (X - 1728)^2 (X^2 - 44031499226496X - 292143758886942437376)\\
&= (X-1728)^2 H_{-100}(X),
\end{align*}
by Lemma 2.  The roots of $H_{-100}(X)$ are supersingular (for $p \ge 7$) exactly when $\left(\frac{-100}{p}\right) = \left(\frac{-4}{p}\right) = -1$, i.e. when $\varepsilon = 1$. Moreover, the factor $Y^2-540Y-6480$ occurs to only the first power in
\begin{align*}
&\textrm{Res}_X(H_{-100}(X),R_5(X,Y)) = (Y^2 - 540Y - 6480)h_{100}(Y)\\
& =(Y^2 - 540Y - 6480) (Y^{10} - 44031499224660Y^9 - 292192545788083696320Y^8\\
& - 111045241276874215905600Y^7 - 64831872214747570823193600Y^6\\
& - 35633053922822368233233495040Y^5 - 19661658654621205173476830924800Y^4\\
& + 2016600435462530152592430229094400Y^3\\
& - 67253379407769529512012174852096000Y^2\\
& + 1082713527360852989716901652332544000Y\\
& - 13177845369975884011784483478416916480),
\end{align*}
for primes not dividing
$$\textrm{Res}_Y( Y^2 - 540Y - 6480,h_{100}(Y)) = 2^{68} \cdot 3^{42} \cdot 5^2 \cdot 7^{12} \cdot 11^6 \cdot 19^4 \cdot 23^2 \cdot 47^2 \cdot 59^2\cdot 71^2 \cdot 83^2.$$
Hence, when $\varepsilon = 1$ and
$$p \notin \{2,3,5,7,11,19,23,47,59,71,83\},$$
the factor $Y^2 - 540Y - 6480$ occurs to exactly the third power in (6): twice for $x = 1728$ and once for $H_{-100}(X)$, when these are supersingular.  This explains the factor $(Y^2 - 540Y - 6480)^\varepsilon$ in (7), since $Y^2 - 540Y - 6480$ exactly divides $ss_p^{(5*)}(Y)$.  \medskip

The multiple roots $y$ corresponding to $x=0$ arise from
$$R_5(0,Y) = (Y^2+216Y+144)^3,$$
while
\begin{align*}
\textrm{Res}_Y&(R_5(X,Y),Y^2+216Y+144) = X^2 H_{-75}(X)\\
& = X^2(X^2 + 654403829760X + 5209253090426880).
\end{align*}
The roots of $H_{-75}(X)$ are supersingular (for $p \ge 7$) exactly when $\left(\frac{-75}{p}\right) = \left(\frac{-3}{p}\right) = -1$, i.e. when $\delta = 1$.  Further,
\begin{align*}
&\textrm{Res}_X(H_{-75}(X),R_5(X,Y)) = (Y^2 +216Y +144)h_{75}(Y)\\
& =(Y^2 +216Y +144) (Y^{10} + 654403830840Y^9 + 5439603238969680Y^8\\
& - 1949338201633113600Y^7 + 473463907652088230400Y^6\\
& - 104049869016988552310784Y^5 + 22874519246403909048606720Y^4\\
& - 1239769622718575548557557760Y^3 + 51906281918763496455571046400Y^2\\
& - 194668556748911160362178969600Y + 279141650822621456977854726144),
\end{align*}
where
$$\textrm{Res}_Y(Y^2 +216Y +144,h_{75}(Y)) = 2^{102} \cdot 3^{26} \cdot 5^2 \cdot 11^8 \cdot 17^2 \cdot 23^2 \cdot 47^2 \cdot 59 \cdot 71.$$
It follows that the exact power of $Y^2 +216Y +144$ dividing (6) is the fourth, when $\delta = 1$ and $p >71$, which explains the presence of the factor $(Y^2 +216Y +144)^{2\delta}$ in (7).  \bigskip

Facts 1-4 prove the equality in (6) and (7) for all primes $p$ not in the set
$$S_5 = \{2,3,5,7,11,13,17,19, 23, 43, 47, 59, 67, 71, 83\}.$$
Using Tables 1 and 2 we check Theorem 1 directly for the 12 primes $\ge 7$ in $S_5$.  This completes the proof of Theorem 1. \bigskip

\begin{table}
  \centering 
  \caption{$ss_p(x)$ for $3 < p < 100$.}\label{tab:1}

\noindent \begin{tabular}{|c|l|c|}
\hline
   &   \\
$p$	&   $ss_p(x) \ \textrm{mod} \ p$  \\
\hline
& \\
5 & $x$\\
7 & $x+1$\\
11 & $x(x+10)$ \\
13 & $x+8$ \\
 17  &  $x(x+9)$  \\
19  &  $(x+1)(x+12)$  \\
23 & $x(x+4)(x+20)$ \\
29 & $x(x+4)(x+27)$ \\
31 & $(x+8)(x+27)(x+29)$ \\
37 & $(x + 29)(x^2 + 31x + 31)$\\
41 & $x(x+9)(x+13)(x+38)$ \\
 43 & $(x+2)(x+35)(x^2+19x+16)$ \\
47 & $x(x+3)(x+11)(x+37)(x+38)$ \\
53 & $x(x + 3)(x + 7)(x^2 + 50x + 39)$ \\
 59 & $x(x + 11)(x + 12)(x + 31)(x+42)(x + 44)$ \\
 61 & $(x + 11)(x + 20)(x + 52)(x^2 + 38x + 24)$\\
 67 & $(x + 1)(x + 14)(x^2 + 8x + 45)(x^2 + 44x + 24)$\\
 71 & $x(x + 5)(x + 23)(x + 30)(x + 31)(x + 47)(x + 54)$ \\
 73 & $(x+17)(x+64)(x^2+57x+8)(x^2+68x+9)$\\
 79 & $(x + 10)(x + 15)(x + 58)(x + 62)(x + 64)(x^2 + 14x + 62)$\\
 83 & $x(x + 15)(x + 16)(x + 33)(x + 55)(x + 66)(x^2 + 7x + 73)$ \\
 89 & $x(x + 23)(x + 37)(x + 76)(x + 82)(x + 83)(x^2 + 26x + 56)$ \\
 97 & $(x+77)(x+96)(x^2+7x+45)(x^2+32x+67)(x^2+42x+8)$\\
 & \\
  \hline
\end{tabular}
\end{table}

\begin{table}
  \centering 
  \caption{$ss_p^{(5*)}(Y)$ for $p \in S_5 - \{2, 3, 5\}$.}\label{tab:1}

\noindent \begin{tabular}{|c|l|c|}
\hline
   &   \\
$p$	&   $ss_p^{(5*)}(Y) \ \textrm{mod} \ p$  \\
\hline
& \\
7 & $Y(Y+3)$\\
11 & $Y(Y+3)(Y+4)(Y+7)$ \\
13 & $(Y + 3)(Y + 9)(Y^2 + 8Y + 10)$ \\
 17  &  $(Y + 13)(Y^2 + 7Y + 1)(Y^2 + 12Y + 8)$  \\
19  &  $Y(Y+1)(Y+2)(Y+9)(Y+11)(Y+15)$  \\
23 & $Y(Y + 5)(Y^2 + 9Y + 6)(Y^2 + 12Y + 6)$ \\
 43 & $Y(Y + 3)(Y + 14)(Y + 25)(Y + 28)(Y + 39)(Y + 41)$\\
 & \ \ \ $\times (Y^2 + 6Y + 40)(Y^2 + 19Y + 13)$ \\
47 & $Y(Y + 29)(Y^2 + 12Y + 3)(Y^2 + 17Y + 2)$\\
& \ \ \ $\times (Y^2 + 24Y + 6)(Y^2 + 28Y + 3)(Y^2 + 34Y + 2)$ \\
 59 & $Y(Y + 3)(Y + 16)(Y + 19)(Y + 20)(Y + 23)(Y + 28)(Y + 30)(Y + 41)$\\
  & \ \ \ $\times (Y + 58)(Y^2 + 15Y + 1)(Y^2 + 24Y + 35)(Y^2 + 58Y + 51)$ \\
 67 & $Y(Y + 3)(Y + 12)(Y + 25)(Y + 28)(Y + 31)(Y + 49)(Y + 54)(Y + 62)$\\
 & \ \ \ $\times (Y^2 + 14Y + 47)(Y^2 + 20Y + 47)(Y^2 + 44Y + 16)(Y^2 + 63Y + 19)$\\
 71 & $Y(Y + 1)(Y + 2)(Y + 3)(Y + 6)(Y + 21)(Y + 26)(Y + 35)(Y + 53)(Y + 66)$\\
 & \ \ \ $\times (Y + 70)(Y^2 + 3Y + 6)(Y^2 + 11Y + 9)(Y^2 + 12Y + 2)(Y^2 + 27Y + 27)$ \\
 83 & $Y(Y + 11)(Y + 39)(Y + 65)(Y + 79)(Y^2 + 12Y + 31)(Y^2 + 23Y + 28)$\\
 & \ \ \ $\times (Y^2 + 24Y + 1)(Y^2 + 35Y + 26)(Y^2 + 41Y + 77)$\\
 & \ \ \ $\times (Y^2 + 50Y + 61)(Y^2 + 57Y + 10)(Y^2 + 65Y + 26)$ \\
 & \\
  \hline
\end{tabular}
\end{table}

\section{The case $N=7$.}

Let the curve $R_7(X,Y) = 0$ be defined by
\begin{align*}
R_7(X,Y) &= X^2-XY(Y^2-21Y+8)(Y^4-42Y^3+454Y^2-1008Y-1280)\\
& \ \ +Y^2(Y^2+224Y+448)^3.
\end{align*}
This is a curve of genus $0$ parametrized by the equations
$$X =\frac{(z^2-3z+9)(z^2-11z+25)^3}{z-8}, \ \ \ Y = \frac{z^2-3z+9}{z-8}$$
See \cite[p. 264]{mor1}.  We have
\begin{align}
\notag \textrm{disc}_X R_7(X,Y) &= (Y + 1)(Y - 27)Y^2(Y - 2)^2(Y - 8)^2(Y - 24)^2\\
& \ \ \ \ \ \ \times (Y^2 - 16Y - 8)^2\\
\notag \textrm{disc}_Y R_7(X,Y) & = -7^7 X^6(X - 1728)^4(X - 54000)^2(X + 96^3)^2\\
 & \ \ \ \ \ \ \times (X + 12288000)^2(X^2 - 4834944X + 14670139392)^2\\
\notag & = -7^7 X^6(X-1728)^4 H_{-12}(X)^2 H_{-19}(X)^2 H_{-27}(X)^2 H_{-24}(X)^2.
\end{align}
Define
\begin{equation}
\mu_7 = \frac{1}{2}\left(1-\left(\frac{-7}{p}\right)\right).
\end{equation}

We want to prove the following. \bigskip

\noindent {\bf Theorem 3.} {\it For a prime $p \ge 5$ and $p \neq 7$ we have the following congruence modulo $p$:}
\begin{align}
&(Y+1)^{\mu_7}(Y-27)^{\mu_7} \textrm{Res}_X(ss_p(X),R_7(X,Y)) \equiv \\
& (Y^2+224Y+448)^{2\delta}(Y^4-528Y^3-9024Y^2-5120Y-1728)^\varepsilon ss_p^{(7*)}(Y)^2.
\end{align}

\noindent {\bf Lemma 4.} {\it We have the following class equations:}
\begin{align*}
H_{-7}(X) &= X+15^3,\\
H_{-28}(X) & = X-255^3,\\
H_{-24}(X) & = X^2 - 4834944X + 14670139392,\\
H_{-147}(X) & = X^2 + 34848505552896000X + 11356800389480448000000,\\
H_{-196}(X) & = X^4 - 12626092121367165696X^3\\
& \ \ \ \ - 44864481851299856707307347968X^2\\
& \ \ \ \ +250850701957837760512539510177792 X\\
& \ \ \ \ -2108010653658430719613224868701536256.
\end{align*}

\noindent {\it Proof.} See Cox \cite[p. 237]{co} for $H_{-7}(X), H_{-28}(X)$.  For $H_{-24}(X)$ see Fricke \cite[III, p. 401]{fr} or \cite[p. 1191]{mor}.  One may also use Berwick \cite[p. 57]{berw}, according to which
$$j\left(\sqrt{6}i\right) = 2^6 \cdot 3^3 \cdot (1+\sqrt{2})^5 (-1+3\sqrt{2})^3,$$
and whose minimal polynomial is $H_{-24}(X)$.  From Berwick \cite[pp. 58]{berw} we also have
$$j\left(\frac{-1+7\sqrt{3}i}{2}\right) = -3\sqrt{21} \cdot 2^{15} \cdot15^3 \cdot \left(\frac{5 + \sqrt{21}}{2}\right)^9 (-2 + \sqrt{21})^3,$$
and its minimal polynomial is $H_{-147}(X)$. \medskip

To verify the polynomial $H_{-196}(X)$ we use the same method as in Lemma 2.  The value $r(i)$ has minimal polynomial
$$f(x) = x^4 + 2x^3 - 6x^2 - 2x + 1.$$
This time we use Yi's relation from \cite[Thm. 3.3]{yi} between $u=r(\tau)$ and $v=r(7\tau)$ given by $P_7(u,v)=0$, where
\begin{align*}
P_7(u,v) = & \ u^8 v^7+(-7v^5+1)u^7+7u^6v^3+7(-v^6+v)u^5+35u^4 v^4\\
& + 7(v^7+v^2)u^3-7u^2v^5-(v^8+7v^3)u-v,
\end{align*}
to compute the resultant of $f(t)$ and $P_7(t,y)$:
 \begin{align*}
 &\textrm{Res}_t(f(t), P_7(t,y)) = y^{32} + 6526y^{31} - 560286y^{30} + 1894660y^{29} - 1558920y^{28}\\
&  + 97188y^{27} + 1383158y^{26}- 16089708y^{25} + 33009225y^{24} - 23680900y^{23} \\
&   + 11485610y^{22}+ 17984710y^{21} - 116298560y^{20} + 132435800y^{19}  - 75016500y^{18} \\
& + 109981440y^{17} + 28870465y^{16}- 109981440y^{15} - 75016500y^{14} - 132435800y^{13} \\
& - 116298560y^{12}  - 17984710y^{11} + 11485610y^{10}+ 23680900y^9 + 33009225y^8  \\
& + 16089708y^7+ 1383158y^6 - 97188y^5 - 1558920y^4- 1894660y^3 - 560286y^2 \\
& - 6526y + 1.
 \end{align*}
This is the minimal polynomial $f_{196}(y)$ of $r(7i)$.  Now $H_{-196}(X)$ may be computed using the resultant
\begin{align*}
&\textrm{Res}_y(f_{196}(y),F(y,X)) = 5^{480} (X^4 - 12626092121367165696 X^3\\ 
& - 44864481851299856707307347968X^2 + 250850701957837760512539510177792X\\
& - 2108010653658430719613224868701536256)^8.
\end{align*}
Alternatively, one may use the polynomial $p_{196}(x)$ from \cite[Section 5, Ex. 3]{mor3}:
\begin{align*}
p_{196}(x)=& x^{16} + 14x^{15} + 64x^{14} + 84x^{13} - 35x^{12} - 14x^{11} + 196x^{10} + 672x^9 + 1029x^8\\
& - 672 x^7 + 196x^6 + 14x^5 - 35x^4 - 84x^3 +64x^2 -14x+1;
\end{align*}
which is the minimal polynomial of the value $r\left(\frac{-49+7i}{5}\right)=r\left(\frac{1+7i}{5}\right)$, and compute that
\begin{equation*}
\textrm{Res}_y(p_{196}(y),F(y,X)) = 5^{120}H_{-196}(X)^4.
\end{equation*}
$\square$ \medskip

We turn now to the proof of Theorem 3. \medskip

\noindent {\bf Fact 5.} {\it All the roots of $(Y+1)^{\mu_7}(Y-27)^{\mu_7} \textrm{Res}_X(ss_p(X),R_7(X,Y))$ are roots of $ss_p^{(7*)}(Y)$.}

As in Fact 1 we just have to consider the factor $(Y+1)(Y-27)$ in (13).  We have from Lemma 4 that
$$\textrm{Res}_Y((Y+1)(Y-27),R_7(X,Y)) = (X + 15^3)^2 (X - 225^3)^2= H_{-7}(X)^2 H_{-28}(X)^2.$$
Hence, the factors $Y+1, Y-27$ occur as factors of $ss_p^{(7*)}(Y)$, for $p \neq 7$ if and only if $\left(\frac{-7}{p}\right) = -1$, i.e. 
if and only if $\mu_7 = 1$. \medskip

Furthermore,
\begin{align}
\notag R(-15^3,Y) &= (Y+1)h_{7}(Y) = (Y + 1)(Y^7 + 4046Y^6 - 64799Y^5 + 16442335Y^4\\
& + 14883071Y^3 + 199370017Y^2 - 45950625Y + 11390625),\\
\notag R(225^3,Y) &= (Y-27)h_{28}(Y) = (Y-27)(Y^7 - 16580676Y^6 + 597100245Y^5 \\
\notag & - 6151819849Y^4 + 14341099983Y^3 - 2649367371Y^2 \\
 & - 383438155625Y - 10183036921875).
\end{align}
Since $h_7(-1) = 3^{10} \cdot 5^4 \cdot 7$ and $h_{28}(27) = -3^8 \cdot 5^4 \cdot 7 \cdot 17^4 \cdot 19^2$, then for primes $ p > 19$ the factors $Y+1$ and $Y-27$ occur to exactly the second power in (13) when $\mu_7 = 1$, and so are accounted for by $ss_p^{(7*)}(Y)^2$ in (14). \medskip

\noindent {\bf Fact 6.} {\it Since $R_7(X,Y)$ is quadratic in $X$, each root $y$ of $ss_p^{(7)}(Y) = 0$ arises from exactly two roots $x$ of $R_7(X,y)=0$, except for the following values, which are all roots of the discriminant in equation (10).}

The argument here is similar to the argument in Fact 2:
\begin{align*}
y &= 0 \ \textrm{corresponds to} \ x = 0, \ \textrm{since} \ R_7(0,Y) = Y^2 h_3(Y)^2\\
& \ = Y^2(Y^2 + 224Y + 448)^3 \ \textrm{and} \ R_7(X,0) = X^2;\\
y &= 2 \ \textrm{corresponds to} \ x = 54000, \ \textrm{since} \ R_7(54000,Y) = (Y-2)^2 h_{12}(Y)\\
& = (Y-2)^2(Y^6 - 53324Y^5+ 3340572Y^4 - 47158880Y^3 + 453452848Y^2\\
& \ + 867240000Y + 729000000)\\
& \ \textrm{and} \ R_7(X,2) = (X-54000)^2;\\
y &= 8 \ \textrm{corresponds to} \ x = -96^3, \ \textrm{since} \ R_7(-96^3,Y) = (Y-8)^2 h_{19}^*(Y)\\
& = (Y-8)^2(Y^6 + 885424Y^5- 41419776Y^4 + 481543168Y^3 + 799436800Y^2\\
& \ + 2916089856Y + 12230590464)\\
& \ \textrm{and} \ R_7(X,8) = (X+96^3)^2;\\
y &= 24 \ \textrm{corresponds to} \ x = -12288000, \ \textrm{since} \ R_7(-12288000,Y)\\
& = (Y-24)^2 h_{27}(Y)\\
& = (Y-24)^2(Y^6+ 12288720Y^5 - 184134144Y^4 + 610171904Y^3\\
& \ + 1748692992Y^2 + 21626880000Y + 262144000000)\\
& \ \textrm{and} \ R_7(X,24) = (X+12288000)^2.
\end{align*}
It follows that for these values $(Y-y)^2$ exactly divides the resultant in (13), when the corresponding $X$-value is supersingular (corresponding to quadratic discriminants $d = -3, -12, -19, -27$; see (11)), and so are accounted for in (14) by the factor $ss_p^{(7*)}(Y)^2$.  As in Fact 2, this will be true for the primes which do not divide the following values, which are the values of each of the above four cofactors of $(Y-y)^2$ evaluated at $y$:
\begin{align*}
y=0: \ \ & h_3(0) = 448 = 2^6 \cdot 7;\\
y= 2: \ \ & h_{12}(2) = 3951763200 = 2^8 \cdot 3^6 \cdot 5^2 \cdot 7 \cdot 11^2;\\
y=8: \ \ & h_{19}(8) = 192631799808 = 2^{22} \cdot 3^8 \cdot 7;\\
y=24: \ \ & h_{27}(24) = 46982810828800 = 2^{22} \cdot 5^2 \cdot 7 \cdot 11^2 \cdot 23^2.
\end{align*}

For the last factor $Y^2-16Y-8$ in (10) we have
\begin{equation*}
\textrm{Res}_Y(R_7(X,Y),Y^2-16Y-8) = (X^2 - 4834944X + 14670139392)^2 = H_{-24}(X)^2
\end{equation*}
and
\begin{align*}
&\textrm{Res}_X(H_{-24}(X),R_7(X,Y)) = (Y^2 - 16Y - 8)^2 h_{24}(X)\\
& = (Y^2 - 16Y - 8)^2 (Y^{12} - 4833568Y^{11} + 11571739408Y^{10} - 2012852637952Y^9\\
& + 15204068799424Y^8 + 493204380225536Y^7 + 11141216141178880Y^6\\
& - 31850426719240192Y^5 + 184900908191444992Y^4 + 1598968808958984192Y^3\\
& + 7770514603029626880Y^2 - 2102123472092135424Y + 3362702965323595776);
\end{align*}
where
$$\textrm{Res}_Y(Y^2-16Y-8),h_{24}(Y)) = 2^{54} \cdot 3^{20} \cdot 7^2 \cdot 13^4 \cdot 17^2 \cdot 19^4 \cdot 23^2.$$
Hence, the factor $Y^2-16Y-8$ is also accounted for in the factorization of (14), for primes $p > 23$.

Note that the above $y$-values are distinct for $p>23$, since
\begin{align*}
\textrm{disc}_Y(Y(Y + 1)(Y - 27)(Y - 2)&(Y - 8)(Y - 24)(Y^2 - 16Y - 8)) = \\
& \ 2^{57} \cdot 3^{32} \cdot 5^8 \cdot 7^2 \cdot 11^2 \cdot 17^4 \cdot 19^2 \cdot 23^2.
\end{align*}

\noindent {\bf Fact 7.} {\it For values $y$ corresponding to two distinct $X$-values, both $X$-values are supersingular when one of them is.}

This follows from the factorization
\begin{align*}
R_7\left(X, \frac{z^2-3z+9}{z-8} \right) & = \left(X-\frac{(z^2 - 3z + 9)(z^2 - 11z + 25)^3}{z - 8}\right)\\
& \ \times \left(X - \frac{(z^2 - 3z + 9)(z^2 + 229z + 505)^3}{(z - 8)^7} \right).
\end{align*}
This is because, with
$$z = \frac{8d^3 - 15d^2 - 9d + 8}{d^3 - 8d^2 + 5d + 1},$$
the quantity
\begin{align}
\notag j_7 & = \frac{(z^2 - 3z + 9)(z^2 + 229z + 505)^3}{(z - 8)^7}\\
& = \frac{(d^2 - d + 1)^3(d^6 - 11d^5 + 30d^4 - 15d^3 - 10d^2 + 5d + 1)^3}{(d^3 - 8d^2 + 5d + 1)(d - 1)^7 d^7}
\end{align}
is the $j$-invariant of the Tate normal form for a point of order $7$:
\begin{equation}
E_7: \ \ Y^2+(1+d-d^2)XY+(d^2-d^3)Y=X^3+(d^2-d^3)X^2;
\end{equation}
and
\begin{align}
\notag j_{7,7} & = \frac{(z^2 - 3z + 9)(z^2 - 11z + 25)^3}{z - 8}\\
& = \frac{(d^2 - d + 1)^3(d^6 + 229d^5 + 270d^4 - 1695d^3 + 1430d^2 - 235d + 1)^3}{d(d - 1)(d^3 - 8d^2 + 5d + 1)^7}
\end{align}
is the $j$-invariant of the isogenous curve
\begin{align}
\notag & E_{7,7}: \ Y^2+(1+d-d^2)XY+7(d^2-d^3)Y= X^3 -d (d-1)(7d+6)X^2\\
\notag & \ -6d(d-1)(d^5 -2d^4-7d^3 + 9d^2 -3d+1)X\\
\notag & \ -d(d-1)(d^9-2d^8-34d^7+153d^6-229d^5+199d^4-111d^3+28d^2-7d+1).
\end{align}
The $j$-invariants in (17) and (19) can be verified using the formulas in \cite[p. 42]{si} (in which the formula for $b_2$ should read $b_2 = a_1^2+4a_2$).  The fact that $E_7$ and $E_{7,7}$ are isogenous can be seen using the method of \cite[Section 5]{mor0}.  Let $\tau$ be the following translation automorphism of the function field $F(x,y)$ defined by the equation (18) for $E_7$:
$$(x,y)^\tau = (x,y) + (0,0) = \left(\frac{d^2(d-1)y}{x^2}, \frac{d^4(d-1)^2(x^2-y)}{x^3}\right).$$
Then $\tau$ has order $7$ and by \cite[Prop. 3.4]{mor0} the fixed field inside $F(x,y)$ of the group $\langle \tau \rangle$ is the field $F(u,v)$, where
\begin{align*}
u &= \sum_{i=0}^6{x^{\tau^i}} = \frac{A(x)}{x^2(d^2-d-x)^2(d^3-d^2-x)^2},\\
v & = \sum_{i=0}^6{y^{\tau^i}} = \frac{B(x)+d(d-1)C(x)y}{x^3(d^2-d-x)^3(d^3-d^2-x)^3}.
\end{align*}
The polynomial $A(x)$ is given by
\begin{align*}
A(x) = & \ x^7+d(d - 1)(d^5 - 2d^4 - 7d^3 + 9d^2 - 3d + 1)x^5\\
& \ -d^3(d - 1)^2(4d^4 - 17d^3 + 12d^2 - 5d + 1)x^4\\
& \ +d^4(d - 1)^3(d^5 - 3d^4 - 4d^2 - 3d - 1)x^3\\
& \ -d^6(d - 1)^4(d + 1)(d^2 - 3d - 3)x^2+d^8(d - 1)^5(d^2 - 3d - 3)x\\
& \ +d^{10}(d - 1)^6.
\end{align*}
The polynomials $B(x)$ and $C(x)$ are given by
\begin{align*}
& B(x) = (x^3 + (d^2 - d)x^2 - (d^5 - 3d^4 + 2d^3)x - d^7 + 2d^6 - d^5)\\
& \times (x^3 -4 (d^3 - d^2)x^2 - (d^7 - 7d^6 + 10d^5 - 3d^4 - d^3)x - 2d^8 + 6d^7 - 6d^6 + 2d^5)\\
& \times (x^3 + (d^3 - 5d^2 + 4d)x^2 + (2d^4 - 3d^3 + d)x - d^6 + 3d^5 - 3d^4 + d^3);
\end{align*}
and
\begin{align*}
C(x) &= (d^3 + d - 1)x^9 + (d^7 - 3d^6 - 8d^5 + 13d^4 - 5d^3 + 2d^2 + 2d - 1)x^8\\
& -d^2(d - 1)(6d^6 - 32d^5 + 28d^4 - 15d^3 + 5d^2 + 18d - 2)x^7\\
& + d^3(d - 1)^2(3d^7 - 13d^6 + 4d^5 - 23d^4 - 7d^3 + 52d^2 + 9d + 3)x^6\\
& -d^4(d - 1)^3(d^8 - 5d^7 + 11d^6 - 28d^5 - 44d^4 + 63d^3 + 41d^2 + 16d + 1)x^5\\
& + d^6(d - 1)^4(d^6 + 5d^5 - 52d^4 + 15d^3 + 60d^2 + 36d + 5)x^4\\
& -d^8(d - 1)^5(3d^5 - 10d^4 - 21d^3 + 33d^2 + 41d + 10)x^3\\
&  +d^{10}(d - 1)^6(d^4 - 8d^3 + 2d^2 + 23d + 10)x^2\\
&  + d^{12}(d - 1)^7(2d^2 - 5d - 5)x\\
&  + d^{14}(d - 1)^8.
\end{align*}
A calculation on Maple shows that if $P=(x,y)$ is a point on $E_7$, then $\varphi(P) = (u,v)$ is a point on $E_{7,7}$.  This shows that $\varphi: E_7 \rightarrow E_{7,7}$ is an isogeny, and therefore that $j_7$ is supersingular if and only if $j_{7,7}$ is supersingular. \medskip

\noindent {\bf Fact 8.} {\it The only roots $y$ of $ss_p^{(7*)}(Y)$ which occur to a power higher than the second in (13) are those which correspond to the roots of the discriminant (11).}

We may restrict our attention to the values of $y$ corresponding to $x=0$ and $x=1728$, since the roots $x = 54000, -96^3, 12288000$
 and the roots of $H_{-24}(X)$ have been handled in Fact 6.  As in the discussion of Fact 4 above, the polynomials $h_{12}(Y), h_{19}^*(Y), h_{27}(Y)$ and $h_{24}(Y)$ occur to the first power in the calculations in Fact 6 and have distinct roots for primes not in the set
 \begin{align*}
 \{2,&3,5,7,11,13,17,19,23,29,37,41,43,47,53,\\
 &61,67,71,89,109,113, 137, 139,157,163\}.
 \end{align*}

For $x=0$ we have $R_7(0,Y) = Y^2(Y^2 + 224Y + 448)^3$ and
\begin{align*}
\textrm{Res}_Y(&R_7(X,Y),Y^2 + 224Y + 448) = X^2 H_{-147}(X)\\
&\ = X^2(X^2 + 34848505552896000X + 11356800389480448000000).
\end{align*}
Hence, the factor $h_3(Y)=Y^2 + 224Y + 448$ occurs in $ss_p^{(7*)}(Y)$ if and only if $\left(\frac{-3}{p}\right) = \left(\frac{-147}{p}\right) = -1$, i.e., if and only if $\delta = 1$.  Furthermore,
$$\textrm{Res}_X(H_{-147}(X), R_7(X,Y)) = (Y^2 + 224Y + 448) h_{147}(Y),$$
for a factor $h_{147}(Y)$ of degree $14$ for which
\begin{equation}
\textrm{Res}_Y(g(Y),h_{147}(Y)) = 2^{108} \cdot 3^{32} \cdot 5^{20} \cdot 7^2 \cdot 11^6 \cdot 17^7 \cdot 23^2 \cdot 29^2 \cdot 47 \cdot 71^2 \cdot 83 \cdot 131.
\end{equation}
When $h_3(Y)$ occurs, it occurs to the fourth power: three times for $x=0$ and once for the roots of $H_{-147}(X)$.  This accounts for the factor $(Y^2 + 224Y + 448)^{2\delta}$ in (14), for the primes not dividing (20). \medskip

For $x = 1728$ we note that
$$R_7(1728,Y) = (Y^4 - 528Y^3 - 9024Y^2 - 5120Y - 1728)^2$$
and
\begin{align*}
&\textrm{Res}_Y(R_7(X,Y),Y^4 - 528Y^3 - 9024Y^2 - 5120Y - 1728) = (X - 1728)^4\\
& \ \times (X^4 - 12626092121367165696X^3 - 44864481851299856707307347968X^2\\
& \ \ \ + 250850701957837760512539510177792X\\
& \ \ \ - 2108010653658430719613224868701536256)\\
&  = (X-1728)^4 H_{-196}(X),
\end{align*}
by Lemma 4.  Thus, the factor
$$g(Y) = Y^4 - 528Y^3 - 9024Y^2 - 5120Y - 1728$$
occurs as a factor in (13) if and only if $\left(\frac{-4}{p}\right)=\left(\frac{-196}{p}\right) = -1$, i.e., if and only if $\varepsilon = 1$.  When it occurs, it does so to the third power: twice for $x=1728$ and once for the roots of $H_{-196}(X)$, since
$$\textrm{Res}_X(H_{-196}(X),R_7(X,Y)) = g(Y) h_{196}(Y),$$
for a factor $h_{196}(Y)$ of degree $28$, for which
\begin{align}
\notag \textrm{Res}_Y(g(Y),&h_{196}(Y)) = 2^{276} \cdot 3^{182} \cdot 7^4 \cdot 11^{30} \cdot 19^{14} \cdot 23^{22} \cdot 31^6 \cdot 43^2 \cdot 47^4\\
& \cdot 59^2 \cdot 71^4 \cdot 79^2 \cdot 83^2 \cdot 107^2 \cdot 131^4 \cdot 151^2 \cdot 167^2 \cdot 179^2 \cdot 191^2.
\end{align}
This accounts for the factor $g(Y)^\varepsilon$ in (14), for the primes not dividing the resultant in (21).  \medskip

Taken together, Facts 5-8 prove Theorem 3, for the primes $p$ not in the set
\begin{align*}
S_7 =  \{2, &3, 5, 7, 11, 13, 17, 19, 23, 29, 31, 37,41,43,47,53,59, 61,67,71,79,\\
 &83,89,107, 109,113, 131, 137, 139,151, 157,163, 167, 179, 191\}.
 \end{align*}
 For the $32$ primes $p \in S_7-\{2,3,7\}$ we can check the assertion of Theorem 3 directly.  Table 3 contains the polynomials $ss_p^{(7*)}(Y)$ for the $19$ primes in $S_7-\{2,3,7\}$ satisfying $p \le 83$.  For larger primes $ss_p^{(7*)}(Y)$ can be calculated using the fact that
 \begin{align*}
 ss_p(X) &\equiv X^\delta (X-1728)^\varepsilon J_p(X),\\
J_p(X) &\equiv \sum_{k=0}^{n_p}{{2n_p + \varepsilon \atopwithdelims ( ) 2k + \varepsilon}{2n_p-2k \atopwithdelims ( ) n_p-k}(-432)^{n_p-k}(t-1728)^k} \ \ (\textrm{mod} \ p),
 \end{align*}
where $n_p = [p/12]$.  (See \cite{mor0}.)  To verify the congruence of Theorem 3 for $p$, it is only necessary to check that the factors which occur to the first power in
$$\textrm{Res}_X(ss_p(X),R_7(X,Y))$$
or to a power higher than the second agree with the extra factors in (13) and (14).  This completes the proof of Theorem 3.  \medskip

\noindent {\bf Corollary.} {\it The degree of $ss_p^{(7*)}(Y)$ is given by}
$$\textrm{deg}(ss_p^{(7*)}(Y)) = \frac{1}{3}\left(p - \left(\frac{-3}{p}\right)\right)+\mu_7.$$

\noindent {\it Proof.} Let $d_p = \textrm{deg}(ss_p^{(7*)}(Y))$.  The formula of Theorem 3 gives directly on taking degrees that
$$2d_p+4\delta+4\varepsilon = 2\mu_7+8 \textrm{deg}(ss_p(X)),$$
since $R_7(X,Y)$ is monic and has degree $8$ in $Y$.  Thus
$$d_p = 4 \textrm{deg}(ss_p(X))-2\delta-2\varepsilon+\mu_7.$$
Now use the fact that
\begin{equation*}
\textrm{deg}(ss_p(X)) = \frac{1}{12}(p-1-4\delta-6\varepsilon)+\delta+\varepsilon.
\end{equation*}
This yields
\begin{align*}
d_p &=\frac{1}{3}(p-1-4\delta-6\varepsilon)+2\delta+2\varepsilon+\mu_7\\
&= \frac{1}{3}(p-1+2\delta)+\mu_7,
\end{align*}
which agrees with the assertion. $\square$ \medskip

The statement in the above corollary is contained in Nakaya's Conjectures 1 and 6 in \cite{na}.

\begin{table}
  \centering 
  \caption{$ss_p^{(7*)}(Y)$ for $p \in S_7 - \{2, 3, 7\}$ and $p \le 83$.}\label{tab:1}

\noindent \begin{tabular}{|c|l|c|}
\hline
   &   \\
$p$	&   $ss_p^{(7*)}(Y) \ \textrm{mod} \ p$  \\
\hline
& \\
5 & $Y(Y+1)(Y+3)$\\
11 & $Y(Y + 9)(Y^2 + 4Y + 8)$ \\
13 & $(Y + 1)(Y + 5)(Y + 12)(Y^2 + 10Y + 5)$ \\
 17  &  $Y(Y + 1)(Y + 7)(Y + 10)(Y + 11)(Y + 13)(Y + 15)$  \\
19  &  $(Y + 1)(Y + 8)(Y + 11)(Y^2 + 3Y + 11)(Y^2 + 4Y + 8)$  \\
23 & $Y(Y + 8)(Y + 21)(Y + 22)(Y^2 + 3Y + 20)(Y^2 + 17Y + 11)$ \\
29 & $Y(Y + 5)(Y + 21)(Y + 27)(Y^2 + 18Y + 8)(Y^2 + 21Y + 13)(Y^2 + 26Y + 12)$ \\
31 & $(Y + 1)(Y + 4)(Y + 8)(Y + 23)(Y + 30)(Y^2 + 4Y + 8)(Y^2 + 20Y + 4)$ \\
& \ \ \ $\times (Y^2 + 23Y + 30)$\\
37 & $(Y + 8)(Y + 14)(Y + 27)(Y + 29)(Y^2 + 21Y + 29)(Y^2 + 23Y + 26)$\\
& \ \ \ $\times (Y^2 + 31Y + 29)(Y^2 + 34Y + 8)$\\
41 & $Y(Y + 1)(Y + 8)(Y + 12)(Y + 13)(Y + 14)(Y + 17)(Y + 29)(Y + 31)$\\
& \ \ \ $\times (Y + 33)(Y + 39)(Y^2 + Y + 18)(Y^2 + 37Y + 26)$ \\
 43 & $(Y + 8)(Y + 27)(Y^2 + 3Y + 8)(Y^2 + 17Y + 41)(Y^2 + 18Y + 42)$\\
 & \ \ \ $\times (Y^2 + 27Y + 35)(Y^2 + 34Y + 11)(Y^2 + 40Y + 11)$ \\
47 & $Y(Y + 1)(Y + 10)(Y + 16)(Y + 20)(Y + 23)(Y + 26)(Y + 31)(Y + 34)$\\
& \ \ \ $\times (Y + 44)(Y + 45)(Y^2 + 15Y + 42)(Y^2 + 26Y + 15)(Y^2 + 27Y + 33)$ \\
53 & $Y(Y + 8)(Y + 9)(Y + 18)(Y + 29)(Y + 45)(Y + 48)(Y + 51)(Y^2 + 23)$\\
& \ \ \ $\times (Y^2 + 12Y + 24)(Y^2 + 13Y + 8)(Y^2 + 37Y + 25)(Y^2 + 50Y + 3)$ \\
 59 & $Y(Y + 1)(Y + 8)(Y + 32)(Y + 35)(Y + 47)(Y + 51)(Y + 52)(Y + 54)(Y + 55)$\\
 & \ \ \ $\times (Y + 57)(Y^2 + 4Y + 8)(Y^2 + 19Y + 23)(Y^2 + 26Y + 14)$\\
  & \ \ \ $\times (Y^2 + 39Y + 50)(Y^2 + 40Y + 40)$ \\
  61 & $(Y + 1)(Y + 3)(Y + 8)(Y + 34)(Y + 58)(Y^2 + 5Y + 9)(Y^2 + 14Y + 38)$\\
  & \ \ \ $\times (Y^2 + 23Y + 58)(Y^2 + 27Y + 53)(Y^2 + 30Y + 34)(Y^2 + 45Y + 53)$\\
  & \ \ \ $\times (Y^2 + 53Y + 33)(Y^2 + 54Y + 28)$\\
  67 & $(Y+8)(Y+59)(Y+62)(Y+64)(Y^2+9 Y+3)(Y^2+27 Y+8)(Y^2+29 Y+45)$\\
  & \ \ \ $\times (Y^2+44 Y+40) (Y^2+51Y+59)(Y^2+58Y+9)(Y^2+62Y+58)$\\
  & \ \ \ $\times (Y^2+66 Y+27) (Y^2+66 Y+52)$\\
 71 & $Y(Y + 47)(Y + 62)(Y + 63)(Y + 64)(Y + 69)(Y^2 + 18)(Y^2 + 4Y + 8)$\\
 & \ \ \ $\times (Y^2 + 9Y + 65)(Y^2 + 11Y + 22)(Y^2 + 23Y + 37)(Y^2 + 26Y + 37)$\\
 & \ \ \ $\times (Y^2 + 27Y + 62)(Y^2 + 31Y + 4)(Y^2 + 63Y + 3)$ \\
 79 & $(Y + 62)(Y + 71)(Y^2 + 4Y + 8)(Y^2 + 11Y + 21)(Y^2 + 12Y + 57)(Y^2 + 17Y + 10)$\\
 & \ \ \ $\times (Y^2 + 19Y + 62)(Y^2 + 23Y + 58)(Y^2 + 27Y + 52)(Y^2 + 47Y + 69)$\\
 & \ \ \ $\times (Y^2 + 56Y + 38)(Y^2 + 57Y + 78)(Y^2 + 71Y + 58)(Y^2 + 78Y + 14)$ \\
 83 & $Y(Y + 1)(Y + 3)(Y + 17)(Y + 24)(Y + 34)(Y + 41)(Y + 54)(Y + 56)(Y + 59)$\\
 & \ \ \ $\times (Y + 72)(Y + 74)(Y + 81)(Y^2 + 9Y + 52)(Y^2 + 21Y + 60)$\\
 & \ \ \ $\times (Y^2 + 25Y + 34)(Y^2 + 26Y + 1)(Y^2 + 31Y + 41)(Y^2 + 45Y + 65)$\\
 & \ \ \ $\times (Y^2 + 72Y + 52)(Y^2 + 74Y + 7)$ \\
 & \\
  \hline
\end{tabular}
\end{table}

\section{The cases $N=2$ and $N=3$.}

Let the polynomial $R_2(X,Y)$ be defined by
$$R_2(X,Y) = X^2-X(Y^2-207Y+3456)+(Y+144)^3,$$
where
\begin{align}
\textrm{disc}_X R_2(X,Y) &= Y(Y - 256)(Y - 81)^2\\
\textrm{disc}_Y R_2(X,Y) & = 4X^2(X - 1728)(X + 15^3)^2 = 4X^2(X-1728)H_{-7}(X)^2.
\end{align}
The curve $R_2(X,Y) = 0$ is parametrized by
$$X = \frac{2^8(z^2-z+1)^3}{z^2(z-1)^2}, \ \ \ Y = \frac{16(z+1)^4}{z(z-1)^2}.$$

Similarly, the polynomial
$$R_3(X,Y) = X^2 - XY(Y^2-126Y+2944)+Y(Y+192)^3,$$
has
\begin{align}
\textrm{disc}_X R_3(X,Y) &= Y(Y - 108)(Y - 8)^2(Y - 64)^2,\\ 
\notag \textrm{disc}_Y R_3(X,Y) & = -27X^2(X - 1728)^2(X - 8000)^2(X + 32768)^2\\
& = -27X^2(X-1728)^2 H_{-8}(X)^2 H_{-11}(X)^2;
\end{align}
and the curve $R_3(X,Y)= 0$ is parametrized by
$$X = \frac{z^3(z^3-24)^3}{z^3-27}, \ \ \ Y = \frac{z^6}{z^3-27}.$$
Also, set
$$\mu_2 = \frac{1}{2}\left(1-\left(\frac{-2}{p}\right)\right).$$ \smallskip

\noindent {\bf Theorem 5.} {\it The following formulas hold for primes $p \ge 5$}:
\begin{align}
\notag Y^\varepsilon (Y-256)^{\mu_2} \textrm{Res}_X&(ss_p(X),R_2(X,Y))\\
& \equiv (Y+144)^{2\delta} (Y-648)^\varepsilon ss_p^{(2*)}(Y)^2 \ \ (\textrm{mod} \ p);\\
\notag & \\
\notag Y^\delta (Y-108)^\delta \textrm{Res}_X&(ss_p(X),R_3(X,Y))\\
& \equiv (Y+192)^{2\delta} (Y^2-576Y-1728)^\varepsilon ss_p^{(3*)}(Y)^2 \ \ (\textrm{mod} \ p).
\end{align}

\noindent {\it Proof of (26).} Formula (26) is proved according to the pattern established for the proofs of Theorems 1 and 3. \medskip

\noindent 1. The roots of the left side of (26) are roots of $ss_p^{(2*)}(X)$ when $\varepsilon = 1$, respectively $\mu_2 = 1$, since $R_2(1728,0) = 0$ and $1728$ is supersingular when $\varepsilon = 1$; and $R_2(20^3,256) = 0$, where $20^3$ is supersingular when $\mu_2 = 1$, since $H_{-8}(X) = X - 20^3$.  (See Cox, \cite[p. 23]{co}.) \medskip

\noindent 2. The values of $Y$ arising from only one value of $X$ are the roots of (22):
\begin{align*}
y &= 0 \ \textrm{corresponds to} \ x = 1728, \ \textrm{since} \ R_2(1728,Y) = Y(Y-648)^2\\
& \ \textrm{and} \ R_2(X,0) = (X-1728)^2;\\
y &= 256 \ \textrm{corresponds to} \ x = 20^3, \ \textrm{since} \ R_2(20^3,Y) = (Y - 256)h_8(Y)\\
& = (Y-256)(Y^2 - 7312Y - 153664) \ \textrm{and} \ R_2(X,256) = (X-20^3)^2;\\
y &= 81 \ \textrm{corresponds to} \ x = -15^3, \ \textrm{since} \ R_2(-15^3,Y) = (Y - 81)^2(Y + 3969)\\
&  \ \textrm{and} \ R_2(X,81) = (X+15^3)^2.
\end{align*}
All other roots of the left side of (26) occur for two distinct values of $x$.  Note that $Y$ and $Y-256$ occur to exactly the first power in the resultant in (26), when $p \notin \{2, 3, 5, 7\}$, since $0$ and $256$ are not roots of the respective cofactors for these primes. This explains the factors $Y^\varepsilon$ and $(Y-256)^{\mu_2}$ in (26). \medskip

\noindent 3. The roots of
$$R_2\left(X,\frac{16(z+1)^4}{z(z-1)^2}\right) = \left(X- \frac{2^8(z^2-z+1)^3}{z^2(z-1)^2}\right)\left(X - \frac{16(z^2 + 14z + 1)^3}{z(z - 1)^4}\right)$$
are the $j$-invariants
$$j_2 = j(E_2) = \frac{2^8(z^2-z+1)^3}{z^2(z-1)^2} \ \ \textrm{and} \ \ j_2'=j(E_2') = \frac{16(z^2 + 14z + 1)^3}{z(z - 1)^4}$$
of the respective elliptic curves
\begin{align*}
& E_2: \ Y^2 = X(X-1)(X-1+z),\\
& E_2': \ V^2 = (U - 1 + z)(U^2 - 4U - 4z + 4).
\end{align*}
Furthermore, the formulas
\begin{equation*}
u = \frac{x^2 + z - 1}{x - 1}, \ \ \ v = \frac{(x^2 - 2x - z + 1)y}{(x - 1)^2}
\end{equation*}
define an isogeny from $E_2$ to $E_2'$.  Thus, the values $j_2, j_2'$ are both supersingular when one is.  \medskip

\noindent 4. The factors $Y-y$ which occur to a power higher than the second in (26) correspond to the roots $x$ of (23).  For $x=0$ we have
$R_2(0,Y) = (Y + 144)^3$ and
$$R_2(X,-144) = X(X - 54000) = X H_{-12}(X);$$
where
$$R_2(54000,Y) = (Y + 144)(Y^2 - 53712Y + 18974736).$$
Thus, $Y+144$ occurs to the fourth power when $p$ does not divide
$$\textrm{Res}_Y(Y+144, Y^2 - 53712Y + 18974736) = 2^4 \cdot 3^5 \cdot 5^4 \cdot 11$$
and $\left(\frac{-3}{p}\right) = \left(\frac{-12}{p}\right) = -1$, i.e. $\delta = 1$; this explains the factor $(Y+144)^{2\delta}$ in (26). \medskip

For $x = 1728$ we have $R_2(1728,Y) = Y(Y - 648)^2$ and
\begin{align*}
&R_2(X,648) = (X- 1728)(X-66^3) = (X-1728)H_{-16}(X),\\
&R_2(66^3,Y) = (Y - 648)(Y^2 - 286416Y - 126023688),\\
&\textrm{Res}_Y(Y - 648,Y^2 - 286416Y - 126023688) = -2^3 \cdot 3^8 \cdot 7^2 \cdot 11^2.
\end{align*}
Hence, $Y-648$ occurs to exactly the third power in (26), for primes $p \notin \{2,3,7,11\}$, when $\left(\frac{-4}{p}\right) = \left(\frac{-16}{p}\right) = -1$, i.e., when $\varepsilon = 1$.  This explains the factor $(Y-648)^\varepsilon$ in (26). \medskip

The last root $x = -15^3$ has been handled in 2.  It only remains to check formula for the primes $ p = 5, 7, 11$.  This can be checked directly:
\begin{align*}
(Y-216)\textrm{Res}_X(X,R_2(X,Y)) &\equiv (Y + 4)^4 \equiv (Y+4)^2 ss_5^{(2*)}(X)^2 \ (\textrm{mod} \ 5);\\
Y(Y+3) \textrm{Res}_X(X+1,R_2(X,Y)) &\equiv Y^2(Y + 3)^3 \equiv (Y+3) ss_7^{(2*)}(X)^2 \ (\textrm{mod} \ 7);\\
Y \textrm{Res}_X(X(X+10),R_2(X,Y)) &\equiv Y^2(Y + 1)^5 \equiv (Y+1)^3 ss_{11}^{(2*)}(X)^2 \ (\textrm{mod} \ 11).
\end{align*}
This completes the proof of (26). \medskip

\noindent {\it Proof of (27).}

\noindent 5. The values $y=0$ and $y = 108$ of the left side of (27) are roots of $ss_p^{(3*)}(Y)$ when $\delta = 1$, since
$$R_3(X,0) = X^2 \ \ \textrm{and} \ \ R_3(X,108) = (X - 54000)^2 = H_{-12}(X)^2.$$

\noindent 6. The values of $Y$ arising from only one value of $X$ are the roots of (24):
\begin{align*}
y &= 0 \ \textrm{corresponds to} \ x = 0, \ \textrm{since} \ R_3(0,Y) = Y(Y + 192)^3 \ \textrm{and} \ R_3(X,0) = X^2;\\
y &= 108 \ \textrm{corresponds to} \ x = 54000, \ \textrm{since}\\
& \ \ R_3(54000,Y) = (Y - 108)(Y^3 - 53316Y^2 + 1156464Y - 27000000)\\
& \ \ \textrm{and} \ R_3(X,108) = (X-54000)^2;\\
y &= 8 \ \textrm{corresponds to} \ x = 20^3, \ \textrm{since} \ R_3(20^3,Y) = (Y - 8)^2h_8^*(Y)\\
& \ = (Y - 8)^2(Y^2 - 7408Y + 1000000) \ \textrm{and} \ R_3(X,8) = (X - 20^3)^2;\\
y & = 64  \ \textrm{corresponds to} \ x = -2^{15}, \ \textrm{since} \ R_3(-2^{15},Y) = (Y - 64)^2h_{11}(Y)\\
& \ = (Y - 64)^2(Y^2 + 33472Y + 262144) \ \textrm{and} \ R_3(X,64) = (X +2^{15})^2.
\end{align*}
All other roots of the left side of (27) occur for two distinct values of $x$.  Note that $Y$ and $Y-108$ occur to exactly the first power in the resultant in (27), when $p \notin \{2, 3, 5, 11\}$, since $0$ and $108$ are not roots of the respective cofactors for these primes. This explains the factors $Y^\delta$ and $(Y-108)^\delta$ in (27). \medskip

\noindent 7. The roots of the polynomial
$$R_2\left(X,\frac{z^6}{z^3-27}\right) = \left(X- \frac{z^3(z^3-24)^3}{z^3-27}\right)\left(X - \frac{z^3(z^3+216)^3}{(z^3-27)^3}\right),$$
namely,
$$j_3 = \frac{z^3(z^3-24)^3}{z^3-27} \ \ \textrm{and} \ \ j_3' = \frac{z^3(z^3+216)^3}{(z^3-27)^3},$$
are the $j$-invariants of the isogenous elliptic curves
$$E_3: \ Y^2 + z XY+Y = X^3 \ \ \textrm{and} \ \ E_3': \ V^2 +z UV +3V =U^3 - 6zU - z^3 - 9,$$
by \cite[p. 252]{mor4}.  Thus, the values $j_3, j_3'$ are both supersingular when one is.  \medskip

\noindent 8. The factors $Y-y$ which occur to a power higher than the second in (27) correspond to the roots $x$ of (25).  For $x=0$ we have
$R_3(0,Y) = Y(Y+192)^3$ and
\begin{align*}
 R_3(X,-192) &= X(X + 12288000) = XH_{-27}(X),\\
 R_3(-12288000,Y) &= (Y + 192)(Y^3 + 12288384Y^2\\
& \ - 3907547136Y + 786432000000), \\
 \textrm{Res}_Y(Y + 192,&Y^3 + 12288384Y^2 - 3907547136Y + 786432000000)\\
 & = 2^{22} \cdot 3 \cdot 5^4 \cdot 11 \cdot 23.
\end{align*}
Hence, $Y+192$ occurs to the fourth power in (27) when $p \notin \{2, 3, 5, 11, 23\}$ and $\left(\frac{-3}{p}\right) = \left(\frac{-27}{p}\right) = -1$, i.e. $\delta = 1$; this explains the factor $(Y+192)^{2\delta}$ in (27). \medskip

For $x = 1728$ we have $R_3(1728,Y) = (Y^2 - 576Y - 1728)^2$, where
\begin{align*}
\textrm{Res}_Y&(R_3(X,Y),Y^2 - 576Y - 1728) = (X - 1728)^2H_{-36}(X)\\
& =  (X - 1728)^2(X^2 - 153542016X - 1790957481984),\\
\textrm{Res}_X&(H_{-36}(X),R_3(X,Y)) = (Y^2 - 576Y - 1728) h_{36}(Y)\\
& =(Y^2 - 576Y - 1728)(Y^6 - 153540288Y^5 - 1948490040384Y^4 \\
& - 677563234836480Y^3 - 408250635513974784Y^2\\
& + 53661008686742765568Y - 1856208739742169956352),
\end{align*}
and
$$\textrm{Res}_Y(Y^2 - 576Y - 1728,h_{36}(Y)) = 2^{58} \cdot 3^6 \cdot 7^{12} \cdot 11^6 \cdot 19^2 \cdot 23^2 \cdot 31^2.$$
Now the fact that $X^2 - 153542016X - 1790957481984 = H_{-36}(X)$ follows from \cite[p. 57]{berw} or \cite[p. 201]{fr1}; according to the latter reference,
$$j(3i) = 2^4 \cdot 3\sqrt{3}(1+\sqrt{3})^4(1+2\sqrt{3})^3(2+3\sqrt{3})^3,$$
which is a root of the above quadratic.  It follows that $Y^2 - 576Y - 1728$ divides (27) to the third power, when $\left(\frac{-4}{p}\right) = \left(\frac{-36}{p}\right) = -1$, i.e. $\varepsilon= 1$; this explains the factor $(Y^2 - 576Y - 1728)^\varepsilon$ in (27).
\medskip

The remaining values $x=20^3$ and $-2^{15}$ have been discussed in point 6 above.  The corresponding factors $Y-8$ and $Y-64$ occur to exactly the second power in (27) for primes $p \notin \{2,3,5,7\}$.  This proves (27) for primes $p$ not in the set
$$S_3 = \{2,3,5,7, 11, 19, 23, 31\}.$$
For these primes (27) can be checked directly using the supersingular polynomials in Table 1.  This completes the proof of Theorem 5.

\section{Proof of Nakaya's Conjecture 2.}

\noindent {\bf Theorem 6.} (a) {\it The polynomial $ss_p^{(5*)}(X)$ splits into linear factors over $\mathbb{F}_p$ if and only if $p \in \{2,3,5,7,11,19\}$, i.e., if and only if $p$ divides the order of the Harada-Norton group $HN$.} \smallskip

\noindent (b) {\it The polynomial $ss_p^{(7*)}(X)$ splits into linear factors over $\mathbb{F}_p$ if and only if $p \in \{2,3,5,7,17\}$, i.e., if and only if $p$ divides the order of the Held group $He$.} \medskip

\noindent {\it Proof.} (a) The roots of $ss_p^{(5*)}(X)$ are the roots $y$ of the polynomial
\begin{align*}
&R_5(x,Y) = Y^6 + (-x + 648)Y^5 + (80x + 140400)Y^4 + (-1890x + 10264320)Y^3\\
& + (12600x + 20217600)Y^2 + (-7776x + 13436928)Y + x^2 - 3456x + 2985984,
\end{align*}
as $x$ ranges over the roots of $ss_p(X)$.  If all the roots of $R_5(x,Y)$ lie in $\mathbb{F}_p$, then the coefficients certainly lie in $\mathbb{F}_p$; and considering the coefficient of $Y^5$ shows that $x\in \mathbb{F}_p$, for all supersingular $j$-invariants $x$.  Thus, $p$ can only be one of the primes in the set
$$\mathfrak{S} = \{2, 3, 5, 7, 11, 13, 17, 19, 23, 29, 31, 41, 47, 59, 71\}.$$
Direct computation using Theorem 1 and the polynomials in Table 1 shows that $p$ is one of the $6$ primes in the assertion.  Also see \cite[Table 10]{mor2}.  \medskip

The proof of (b) is the same using
\begin{align*}
R_7(x,Y) = \ & Y^8 + (-x + 672)Y^7 + (63x + 151872)Y^6 + (-1344x + 11841536)Y^5\\
 &+ (10878x + 68038656)Y^4 + (-23520x + 134873088)Y^3\\
 & + (-18816x + 89915392)Y^2 + 10240xY + x^2
\end{align*}
and Theorem 3.
$\square$  \bigskip

The same argument can be used to prove Nakaya's Theorem 5 in \cite{na}, using the fact that the coefficients of $Y^2$ and $Y^3$ in the respective polynomials $R_2(X,Y)$ and $R_3(X,Y)$ are $-X$ plus a constant.  This eliminates the need to use any class number estimates.

\medskip

\noindent Dept. of Mathematical Sciences, LD 270

\noindent Indiana University - Purdue University at Indianapolis (IUPUI)

\noindent Indianapolis, IN 46202

\noindent {\it e-mail: pmorton@iupui.edu}

\end{document}